\documentclass{ifacconf}

\usepackage{graphicx}      % include this line if your document contains figures
\usepackage{natbib}        % required for bibliography

\usepackage{color}

\usepackage{subfigure}

\usepackage{url}
\usepackage{algorithm}
\usepackage{comment}
\usepackage{amsfonts,amssymb,amsmath,mathtools}
\usepackage{latexsym}
\newtheorem{theorem}{\bf Theorem} \newtheorem{definition}[theorem]{\bf Definition} 
\newtheorem{lemma}[theorem]{\bf Lemma} \newtheorem{remark}[theorem]{\bf Remark}
  \newtheorem{proposition}[theorem]{\bf Proposition} 
\newtheorem{assumption}[theorem]{\bf Assumption}  
\newtheorem{Algorithm}[theorem]{\bf Algorithm}

\usepackage{setspace}
\begin{document}
\begin{frontmatter}

\title{On the design of terminal ingredients for data-driven MPC} 

\thanks[footnoteinfo]{This work was funded by Deutsche Forschungsgemeinschaft (DFG, German Research Foundation) under Germany's Excellence Strategy - EXC 2075 - 390740016. The authors thank the International Max Planck Research School for Intelligent Systems (IMPRS-IS) for supporting Julian Berberich, and the International Research Training Group Soft Tissue Robotics  (GRK 2198/1 - 277536708).\\
\textcopyright 2021 the authors. This work has been accepted to IFAC for publication under a 
Creative Commons Licence CC-BY-NC-ND.}

\author[First]{Julian Berberich} 
\author[First]{Johannes K\"ohler} 
\author[Second]{Matthias A. M\"uller}
\author[First]{Frank Allg\"ower}

\address[First]{University of Stuttgart, Institute for Systems Theory and Automatic Control, 70550 Stuttgart, Germany (email:$\{$julian.berberich, johannes.koehler, frank.allgower$\}$@ist.uni-stuttgart.de).}
\address[Second]{Leibniz University Hannover, Institute of Automatic Control, 30167 Hannover, Germany (e-mail:mueller@irt.uni-hannover.de).}

\begin{abstract}
We present a model predictive control (MPC) scheme to control linear time-invariant systems using only measured input-output data and no model knowledge.
The scheme includes a terminal cost and a terminal set constraint on an extended state containing past input-output values.
We provide an explicit design procedure for the corresponding terminal ingredients that only uses measured input-output data.
Further, we prove that the MPC scheme based on these terminal ingredients exponentially stabilizes the desired setpoint in closed loop.
Finally, we illustrate the advantages over existing data-driven MPC approaches with a numerical example.
\end{abstract}

%\begin{keyword}
%Learning and predictive control, output feedback predictive control, stability and recursive feasibility
%\end{keyword}

\end{frontmatter}
%===============================================================================

\section{Introduction}
Designing model predictive control (MPC) schemes based on measured data with only partial or no model knowledge is an active field of research using, e.g., model adaptation~\citep{adetola2011robust,aswani2013provably,tanaskovic2014adaptive} or reinforcement learning based approaches~\citep{berkenkamp2017safe,zanon2020safe}.
A key motivation for such approaches is that obtaining accurate model knowledge can be difficult in practice, whereas large amounts of data are often available and can be exploited for control.
In this paper, we investigate MPC of unknown linear time-invariant (LTI) systems using only measured input-output data and no model knowledge.
Our approach relies on a result by~\cite{willems2005note}, which shows that all trajectories of an LTI system can be parametrized via one persistently exciting data trajectory.
This fact has been used, e.g., for data-driven simulation by~\cite{markovsky2008data} and, subsequently, for data-driven MPC by~\cite{yang2015data,coulson2019deepc}.
While~\cite{coulson2020distributionally} derive open-loop robustness guarantees of the scheme,~\cite{berberich2021guarantees} prove closed-loop stability and robustness properties, even if the measured data are affected by noise.
Further theoretical results on robust constraint satisfaction based on noisy data and on a tracking MPC formulation are derived by~\cite{berberich2020constraints} and~\cite{berberich2020tracking}, respectively.
Notably, the existing works with closed-loop guarantees by~\cite{berberich2021guarantees,berberich2020tracking,berberich2020constraints} require terminal equality constraints, which may result in a small region of attraction and a small robustness margin.

In this paper, we propose a novel data-driven MPC scheme with terminal cost and terminal constraints, thereby improving robustness and increasing the region of attraction of the closed loop, in analogy to results from model-based MPC~(\cite{chen1998quasi,rawlings2020model}).
Similar to other works on data-driven MPC based on~\cite{willems2005note}, we assume that only input-output data of the unknown system are available and thus, we employ an extended state vector involving consecutive input-output measurements.
We prove that the MPC scheme with suitable terminal ingredients involving this extended state exponentially stabilizes the closed loop.
Further, using recent results on data-driven control by~\cite{berberich2020combining}, we show how the terminal ingredients can be constructed using only measured data.
A key benefit of the presented approach is that closed-loop stability guarantees can be given in an end-to-end fashion, without prior system identification steps.
Throughout this paper, we assume that the measured data are not affected by noise, but we conjecture that practical stability using noisy measurements can be proven based on recent robustness results for data-driven MPC in~\cite{berberich2021linearpart2}.

%similar to the terminal equality constrained formulation by~\cite{berberich2021guarantees}.

The work by~\cite{dutta2014certification} is conceptually related to our results since it provides a terminal set constraint for (model-based) MPC with input-output models using an implicit characterization of the maximal invariant set.
Moreover,~\cite{abbas2016robust} propose a linear matrix inequality (LMI) based procedure to construct stabilizing terminal ingredients for linear parameter-varying input-output models.
Compared to these model-based results, the proposed MPC scheme and the computation of the terminal ingredients require only one input-output trajectory of an LTI system and no model knowledge.

The remainder of the paper is structured as follows.
After introducing required preliminaries in Section~\ref{sec:prelim}, we present the data-driven MPC scheme, prove closed-loop stability, and provide a design procedure for terminal ingredients in Section~\ref{sec:MPC}.
The approach is applied to a numerical example in Section~\ref{sec:example}, and the paper is concluded in Section~\ref{sec:conclusion}.

\subsubsection*{Notation}
We denote the set of integers in the interval $[a,b]$ by $\mathbb{I}_{[a,b]}$ and the set of nonnegative integers by $\mathbb{I}_{\geq0}$.
Moreover, we write $\lVert x\rVert_2$ and $\lVert A\rVert_2$ for the (induced) $2$-norm of some vector $x$ and matrix $A$.
For matrices $P=P^\top$, $P_2=P_2^\top$, we write $\lambda_{\min}(P)$ for the minimum eigenvalue of $P$ and we define $\lambda_{\min}(P,P_2)\coloneqq\min\{\lambda_{\min}(P),\lambda_{\min}(P_2)\}$, and similarly for $\lambda_{\max}(P)$ and $\lambda_{\max}(P,P_2)$.
We write $P\succ0$ ($P\succeq0$) if $P$ is positive (semi-) definite, and similarly $P\prec0$ ($P\preceq0$) if $P$ is negative (semi-) definite.
Further, we define $\lVert x\rVert_{P}^2\coloneqq x^\top P x$, we denote matrix entries which can be inferred from symmetry by $\star$, and we write $I$ for an identity matrix of appropriate dimension.
For some generic sequence $\{x_k\}_{k=0}^{N-1}$, we introduce $x_{[a,b]}\coloneqq\begin{bmatrix}x_a^\top&\dots&x_b^\top\end{bmatrix}^\top$ and we abbreviate $x\coloneqq x_{[0,N-1]}$.
Finally, we define the Hankel matrix
\begin{align*}
H_L(x)\coloneqq\begin{bmatrix}x_0&x_1&\dots&x_{N-L}\\
x_1&x_2&\dots&x_{N-L+1}\\
\vdots&\vdots&\vdots&\vdots\\
x_{L-1}&x_L&\dots&x_{N-1}\end{bmatrix}.
\end{align*}

\section{Preliminaries}\label{sec:prelim}
In this section, we introduce preliminaries regarding the data-driven system representation (Section~\ref{subsec:prelim_param}) and the extended state-space description (Section~\ref{subsec:prelim_extended}).
Further, we describe the problem setting in Section~\ref{subsec:prelim_problem}.

\subsection{Data-driven system parametrization}\label{subsec:prelim_param}
We consider discrete-time LTI systems of the form
\begin{align}\label{eq:sys}
\begin{split}
x_{k+1}&=Ax_k+Bu_k,\\
y_k&=Cx_k+Du_k
\end{split}
\end{align}
with state $x_k\in\mathbb{R}^n$, input $u_k\in\mathbb{R}^m$, and output $y_k\in\mathbb{R}^p$, all at time $k\in\mathbb{I}_{\geq0}$, where the matrices $(A,B,C,D)$ are unknown.
Throughout this paper, we assume that $(A,B)$ is controllable and $(A,C)$ is observable.
Controllability is required both for applying the main result of~\cite{willems2005note} and for the presented design of terminal ingredients, whereas observability is w.l.o.g. since unobservable modes do not affect the output cost or constraints.
\begin{definition}\label{def:pe}
We say that $\{u_k\}_{k=0}^{N-1}$ with $u_k\in\mathbb{R}^m$ is persistently exciting of order $L$ if $\mathrm{rank}(H_L(u))=mL$.
\end{definition}
The following main result of~\cite{willems2005note} states that a given input-output trajectory $\{u_k^d,y_k^d\}_{k=0}^{N-1}$ can be used to parametrize all system trajectories if the input $u^d$ is persistently exciting (compare~\cite{berberich2020trajectory} for a description in the state-space framework).
\begin{theorem}\label{thm:hankel}
(\cite{willems2005note})
Suppose $\{u_k^d,y_k^d\}_{k=0}^{N-1}$ is a trajectory of~\eqref{eq:sys}, where $u^d$ is persistently exciting of order $L+n$.
Then, $\{\bar{u}_k,\bar{y}_k\}_{k=0}^{L-1}$ is a trajectory of~\eqref{eq:sys} if and only if there exists $\alpha\in\mathbb{R}^{N-L+1}$ such that
\begin{align}\label{eq:thm_hankel}
\begin{bmatrix}H_L(u^d)\\H_L(y^d)\end{bmatrix}\alpha=\begin{bmatrix}\bar{u}\\\bar{y}\end{bmatrix}.
\end{align}
\end{theorem}
Theorem~\ref{thm:hankel} provides a direct non-parametric system description without identifying a model of the system.
Thus, it can be used to predict system trajectories, which we exploit for the MPC scheme presented in this paper.

\subsection{Extended state representation}\label{subsec:prelim_extended}
We define the lag of the system~\eqref{eq:sys} as follows.
\begin{definition}\label{def:lag}
The lag $\underline{l}$ of~\eqref{eq:sys} is the smallest $l\in\mathbb{I}_{[1,n]}$ such that the following observability matrix has rank $n$:
\begin{align}\label{eq:obsv}
\Phi_l\coloneqq\begin{bmatrix}C^\top&(CA)^\top&\dots&(CA^{l-1})^\top\end{bmatrix}^\top.
\end{align}
\end{definition}
For some integer $l$, we define the extended state $\xi_t$ as
\begin{align}\label{eq:extended_state}
\xi_t\coloneqq\begin{bmatrix}u_{[t-l,t-1]}\\y_{[t-l,t-1]}\end{bmatrix},
\end{align}
where $t\geq l$.
We abbreviate the dimension of $\xi_t$ by $n_{\xi}\coloneqq(m+p)l$ and we denote the discrete-time LTI dynamics corresponding to the extended state by
\begin{align}\label{eq:sys_xi}
\begin{split}
\xi_{k+1}&=\tilde{A}\xi_k+\tilde{B}u_k,\\
y_k&=\tilde{C}\xi_k+\tilde{D}u_k
\end{split}
\end{align}
for suitable matrices $\tilde{A}$, $\tilde{B}$, $\tilde{C}$, $\tilde{D}$.
To be precise, the matrices in~\eqref{eq:sys_xi} take the form shown in~\eqref{eq:sys_xi_detailed}, where the matrices $F_i$, $G_i$, $i\in\mathbb{I}_{[1,l]}$, $D$ are unknown.
\begin{figure*}%[b]%% over both columns
\vspace{2pt}
\begin{align}\label{eq:sys_xi_detailed}
\left[\begin{array}{c}
u_{k-l+1}\\\vdots\\u_{k-1}\\u_k\\\hline y_{k-l+1}\\\vdots\\y_{k-1}\\y_k
\end{array}\right]
=
\left[\begin{array}{cccc|cccc}
0&I&\dots&0&0&\dots&\dots&0\\
\vdots&\ddots&\ddots&\vdots&\vdots&\ddots&\ddots&\vdots\\
0&\ddots&\ddots&I&\vdots&\ddots&\ddots&\vdots\\
0&\dots&\dots&0&0&\dots&\dots&0
\\\hline
0&\dots&\dots&0&0&I&\dots&0\\
\vdots&\ddots&\ddots&\vdots&\vdots&\ddots&\ddots&\vdots\\
0&\dots&\dots&0&0&\dots&\dots&I\\
G_l&\dots&\dots&G_1&
F_l&\dots&\dots&F_1
\end{array}\right]
\left[\begin{array}{c}
u_{k-l}\\\vdots\\u_{k-2}\\u_{k-1}\\\hline y_{k-l}\\\vdots\\y_{k-2}\\y_{k-1}
\end{array}\right]
+
\left[\begin{array}{c}
0\\\vdots\\0\\I\\\hline 0\\\vdots\\0\\D
\end{array}\right]
u_k
\end{align}
\noindent\makebox[\linewidth]{\rule{\textwidth}{0.4pt}}
\end{figure*}
It is straightforward to show that~\eqref{eq:sys} and~\eqref{eq:sys_xi} have an equivalent input-output behavior if $l\geq\underline{l}$ (compare~\cite{goodwin2014adaptive,koch2020provably}).
Throughout this paper, we assume that an upper bound $l$ on the lag $\underline{l}$ is available.
The proposed design method (Section~\ref{subsec:MPC_design}) further requires $pl=n$, which is only valid if $\underline{l}$ is known exactly.
%The extended state $\xi_t$ plays a central role in our MPC approach since we assume that only input-output and no state measurements of~\eqref{eq:sys} are available.
%On the other hand, $\xi_t$ can be constructed using only input-output data which allows us to design terminal ingredients for MPC similar to model-based procedures~\citep{chen1998quasi,rawlings2020model}.
Finally, we define the matrix $T_y$ such that $y_{t-1}=T_y\xi_t$, i.e., $T_y=\begin{bmatrix}0&\dots&0&I\end{bmatrix}$.

\subsection{Problem setting}\label{subsec:prelim_problem}

Throughout this paper, the matrices $A$, $B$, $C$, $D$, or equivalently, $F_i$, $G_i$, $D$, are unknown and one input-output trajectory $\{u_k^d,y_k^d\}_{k=0}^{N-1}$ of~\eqref{eq:sys} is available.
We consider pointwise-in-time constraints on the input and output, i.e., $u_t\in\mathbb{U}\subseteq\mathbb{R}^m$ and $y_t\in\mathbb{Y}\subseteq\mathbb{R}^p$ for $t\in\mathbb{I}_{\geq0}$.
Our objective is stabilization of a given input-output setpoint $(u^s,y^s)$ which is an equilibrium in the following sense.
\begin{definition}\label{def:IO_eq}
We say that $(u^s,y^s)\in\mathbb{R}^{m+p}$ is an equilibrium of~\eqref{eq:sys}, if the sequence $\{\bar{u}_k,\bar{y}_k\}_{k=0}^{\underline{l}}$ with $(\bar{u}_k,\bar{y}_k)=(u^s,y^s)$ for all $k\in\mathbb{I}_{[0,\underline{l}]}$ is a trajectory of~\eqref{eq:sys}.
\end{definition}
We assume $(u^s,y^s)\in\mathrm{int}\left(\mathbb{U}\times \mathbb{Y}\right)$ and we denote the corresponding steady-state of~\eqref{eq:sys_xi} by $\xi^s$.
While assuming knowledge of whether an input-output setpoint $(u^s,y^s)\neq0$ is an equilibrium of the system~\eqref{eq:sys} can be restrictive since $A$, $B$, $C$, $D$ are unknown, this condition can be relaxed by introducing artificial setpoints in the MPC scheme which are optimized online (cf.~\cite{limon2008mpc}).

Our contribution can be summarized as follows.
We propose a data-driven MPC approach to control~\eqref{eq:sys} with closed-loop stability guarantees.
The MPC scheme uses~\eqref{eq:thm_hankel} to predict future trajectories based on one persistently exciting input-output trajectory $\{u_k^d,y_k^d\}_{k=0}^{N-1}$.
In contrast to earlier works, e.g., by~\cite{coulson2019deepc,berberich2021guarantees}, we ensure closed-loop stability via terminal ingredients, i.e., via an appropriate terminal set constraint and a terminal cost.
This follows standard arguments from model-based MPC~\citep{chen1998quasi,rawlings2020model}.
The key difficulty is that, since we do not have access to the state $x_t$ or the matrices $A$, $B$, $C$, $D$ in~\eqref{eq:sys}, the terminal ingredients need to be designed for the extended state $\xi_t$ in~\eqref{eq:sys_xi} using only measured data.
Based on recent results on data-driven control by~\cite{berberich2020combining}, we compute terminal ingredients which can be used to implement the MPC scheme.
Finally, we illustrate the advantages w.r.t. earlier works with a numerical example.

\section{Data-driven MPC with terminal ingredients}\label{sec:MPC}
In this section, we present a data-driven MPC scheme with terminal ingredients and closed-loop stability guarantees.
After defining the MPC scheme in Section~\ref{subsec:MPC_scheme}, we prove the desired closed-loop properties in Section~\ref{subsec:MPC_stab}.
Further, in Section~\ref{subsec:MPC_design}, we show how the corresponding terminal ingredients can be computed using only measured data.

\subsection{Proposed MPC scheme}\label{subsec:MPC_scheme}

\begin{subequations}\label{eq:DD_MPC}
Given an input-output data trajectory $\{u_k^d,y_k^d\}_{k=0}^{N-1}$ and initial conditions $\{u_k,y_k\}_{k=t-l}^{t-1}$, the proposed MPC scheme relies on the following open-loop optimal control problem
\begin{align}\label{eq:DD_MPC_cost}
\underset{\substack{\alpha(t)\\\bar{u}(t),\bar{y}(t)}}{\min}\>\>&\sum_{k=0}^{L-1}\lVert \bar{u}_k(t)-u^s\rVert_R^2+\lVert \bar{y}_k(t)-y^s\rVert_Q^2\\\nonumber
&+\lVert\bar{\xi}_L(t)-\xi^s\rVert_P^2\\\label{eq:DD_MPC_hankel}
\text{s.t.}\quad&\begin{bmatrix}\bar{u}_{[-l,L-1]}(t)\\
\bar{y}_{[-l,L-1]}(t)\end{bmatrix}=\begin{bmatrix}H_{L+l}(u^d)\\H_{L+l}(y^d)\end{bmatrix}\alpha(t),\\\label{eq:DD_MPC_init}
&\begin{bmatrix}\bar{u}_{[-l,-1]}(t)\\\bar{y}_{[-l,-1]}(t)\end{bmatrix}=\begin{bmatrix}u_{[t-l,t-1]}\\y_{[t-l,t-1]}\end{bmatrix},\\\label{eq:DD_MPC_term}
&\bar{\xi}_L(t)\in\Xi_f,\>\bar{\xi}_L(t)=\begin{bmatrix}\bar{u}_{[L-l,L-1]}(t)\\\bar{y}_{[L-l,L-1]}(t)\end{bmatrix},\\\label{eq:DD_MPC_constraints}
&\bar{u}_k(t)\in\mathbb{U},\>\>\bar{y}_k(t)\in\mathbb{Y},\>\>k\in\mathbb{I}_{[0,L-1]}.
\end{align}
\end{subequations}
In the cost~\eqref{eq:DD_MPC_cost}, we penalize the distance of the predicted variables w.r.t. the desired setpoint $(u^s,y^s)$, where $Q\succ0$, $R\succ0$ are design parameters.
The constraint~\eqref{eq:DD_MPC_hankel} parametrizes all input-output trajectories $\bar{u}(t)$, $\bar{y}(t)$ of~\eqref{eq:sys_xi} using Theorem~\ref{thm:hankel}.
Moreover, the predicted trajectory is initialized using past $l$ input-output measurements in~\eqref{eq:DD_MPC_init}, and the scheme includes pointwise-in-time constraints on the predicted input and output in~\eqref{eq:DD_MPC_constraints}.
Finally, the matrix $P\succ0$ and the set $\Xi_f$ characterize a suitable terminal cost $\lVert\bar{\xi}_L(t)-\xi^s\rVert_P^2$ and terminal constraint $\bar{\xi}_L(t)\in\Xi_f$ on the predicted extended state at time $L$.
We assume that these terminal ingredients satisfy the following assumption.
\begin{assumption}\label{ass:term_ingredients}
There exist matrices $P=P^\top\succ0$, $K\in\mathbb{R}^{m\times n_{\xi}}$, and a set $\Xi_f\subseteq\mathbb{U}^l\times\mathbb{Y}^l$ with $\xi^s\in\mathrm{int}(\Xi_f)$ such that for all $\xi\in\Xi_f$, $u=u^s+K(\xi-\xi^s)$, and $y=(\tilde{C}+\tilde{D}K)\xi$, we have
\begin{enumerate}
\item[i)] $\tilde{A}\xi+\tilde{B}u\in\Xi_f$,
\item[ii)] $u\in\mathbb{U}$ and $y\in\mathbb{Y}$,
\item[iii)] the following inequality holds:
\begin{align}\label{eq:ass_term_ingredients}
\lVert (\tilde{A}+\tilde{B}K)\xi\rVert_P^2\leq\lVert\xi\rVert_P^2-\lVert\xi\rVert_{K^\top RK}-\lVert y\rVert_{Q}^2.
\end{align}
\end{enumerate}
\end{assumption}
Assumption~\ref{ass:term_ingredients} is a standard condition in model-based MPC to ensure closed-loop stability, compare~\cite{chen1998quasi,rawlings2020model}.
In Section~\ref{subsec:MPC_design}, we show how a matrix $P$ and a set $\Xi_f$ as in Assumption~\ref{ass:term_ingredients} can be constructed using only measured data and no model knowledge.
If $\mathbb{U}$, $\mathbb{Y}$, and $\Xi_f$ are convex polytopes (ellipsoids), then Problem~\eqref{eq:DD_MPC} is a convex (quadratically constrained) quadratic program which can be solved efficiently.
Throughout this paper, $\bar{u}^*(t)$, $\bar{y}^*(t)$, $\alpha^*(t)$ denote the optimal solution of Problem~\eqref{eq:DD_MPC} at time $t\in\mathbb{I}_{\geq0}$, whereas $u_t$, $y_t$ denote the closed-loop input and output at time $t\in\mathbb{I}_{\geq0}$.
Further, we write $J_L^*(\xi_t)$ for the optimal cost of Problem~\eqref{eq:DD_MPC} with initial condition $\xi_t=\begin{bmatrix}u_{[t-l,t-1]}^\top,y_{[t-l,t-1]}^\top\end{bmatrix}^\top$.
Algorithm~\ref{alg:MPC} summarizes the MPC scheme based on repeatedly solving Problem~\eqref{eq:DD_MPC}.

\begin{algorithm}
\begin{Algorithm}\label{alg:MPC}
\normalfont{\textbf{Data-Driven MPC Scheme}}
\begin{enumerate}
\item At time $t$, take the past $l$ measurements $u_{[t-l,t-1]}$, $y_{[t-l,t-1]}$ and solve~\eqref{eq:DD_MPC}.
\item Apply the input $u_t=\bar{u}_0^*(t)$.
\item Set $t=t+1$ and go back to (1).
\end{enumerate}
\end{Algorithm}
\end{algorithm}

\subsection{Closed-loop guarantees}\label{subsec:MPC_stab}
We make the following common assumption.
\begin{assumption}\label{ass:upper_bound}
The optimal value function is quadratically upper bounded, i.e., there exists $c_u>0$ such that $J_L^*(\xi)\leq c_u\lVert\xi\rVert_2^2$ for all $\xi$ such that Problem~\eqref{eq:DD_MPC} is feasible.
\end{assumption}
Assumption~\ref{ass:upper_bound} is, e.g., satisfied if the sets $\mathbb{U}$, $\mathbb{Y}$, $\Xi_f$ are polytopes~\citep{bemporad2002explicit} or if $\mathbb{U}$ and $\mathbb{Y}$ are compact, compare~\cite[Proposition 2.16]{rawlings2020model}.
More generally, Assumption~\ref{ass:upper_bound} holds with a non-quadratic upper bound, thus leading only to \emph{asymptotic} stability guarantees for the closed loop, if $\mathbb{U}$ is compact~\citep[Proposition 2.16]{rawlings2020model}.

\begin{theorem}\label{thm:stability}
Suppose Assumptions~\ref{ass:term_ingredients} and~\ref{ass:upper_bound} hold and the input $u^d$ is persistently exciting of order $L+l+n$.
If Problem~\eqref{eq:DD_MPC} is feasible at initial time $t=0$, then
\begin{enumerate}
\item[i)] it is feasible at any $t\in\mathbb{I}_{\geq0}$,
\item[ii)] the closed-loop trajectory satisfies the constraints, i.e., $u_t\in\mathbb{U}$ and $y_t\in\mathbb{Y}$ for all $t\in\mathbb{I}_{\geq0}$,
\item[iii)] the equilibrium $\xi^s$ is exponentially stable for the resulting closed loop.
\end{enumerate} 
\end{theorem}
\begin{pf}
W.l.o.g., we assume $u^s=0$, $y^s=0$, $\xi^s=0$.
At time $t+1$, we construct a feasible candidate solution $(\bar{u}'(t+1),\bar{y}'(t+1))$ of~\eqref{eq:DD_MPC} by shifting the previously optimal solution $\bar{u}^*(t)$, $\bar{y}^*(t)$ of~\eqref{eq:DD_MPC} and appending the input $\bar{u}_{L-1}'(t+1)=K\begin{bmatrix}\bar{u}^*_{[L-l,L-1]}(t)\\\bar{y}_{[L-l,L-1]}^*(t)\end{bmatrix}$, with $K$ as in Assumption~\ref{ass:term_ingredients}.
Since $(\bar{u}'(t+1),\bar{y}'(t+1))$ is a trajectory of~\eqref{eq:sys}, there exists $\alpha'(t+1)$ satisfying~\eqref{eq:DD_MPC_hankel} by Theorem~\ref{thm:hankel}, i.e., all constraints of~\eqref{eq:DD_MPC} are satisfied.
Denoting the cost of~\eqref{eq:DD_MPC} for this candidate solution by $J_L'(\xi_{t+1})$, we have
\begin{align}\label{eq:thm_stab_proof_cost_decay}
&J_L^*(\xi_{t+1})-J_L^*(\xi_t)\leq J_L'(\xi_{t+1})-J_L^*(\xi_t)\\\nonumber
=&\lVert\bar{u}_{L-1}'(t+1)\rVert_R^2+\lVert\bar{y}_{L-1}'(t+1)\rVert_Q^2-\lVert\bar{u}_0^*(t)\rVert_R^2-\lVert\bar{y}_0^*(t)\rVert_Q^2\\\nonumber
&+\lVert\bar{\xi}_L'(t+1)\rVert_P^2-\lVert\bar{\xi}_L^*(t)\rVert_P^2\stackrel{\eqref{eq:ass_term_ingredients}}{\leq}-\lVert u_t\rVert_R^2-\lVert y_t\rVert_Q^2.
\end{align}
Since $(A,C)$ is observable, the pair $(\tilde{A},\tilde{C})$ is detectable and hence, there exists an input-output-to-state stability (IOSS) Lyapunov function $W(\xi)=\xi^\top P_W\xi$ satisfying
\begin{align}\label{eq:thm_stab_proof_IOSS_decay}
W(\tilde{A}\xi+\tilde{B}u)-W(\xi)\leq -\frac{1}{2}\lVert\xi\rVert_2^2+c_1\lVert u\rVert_2^2+c_2\lVert y\rVert_2^2
\end{align}
for suitable $c_1,c_2>0$, $P_W\succ0$, and all $u\in\mathbb{R}^m$, $\xi\in\mathbb{R}^{n_\xi}$, $y=\tilde{C}\xi+\tilde{D}u$, compare~\cite{cai2008input}.
We consider as a Lyapunov function candidate the weighted sum of $J_L^*$ and $W$, i.e., $V(\xi_t)=J_L^*(\xi_t)+\gamma W(\xi_t)$ with $\gamma=\frac{\lambda_{\min}(Q,R)}{\max\{c_1,c_2\}}>0$.
Combining~\eqref{eq:thm_stab_proof_cost_decay} and~\eqref{eq:thm_stab_proof_IOSS_decay}, this implies $V(\xi_{t+1})-V(\xi_t)\leq-\frac{\gamma}{2}\lVert\xi_t\rVert_2^2$.
The Lyapunov function candidate $V(\xi_t)$ is clearly quadratically lower bounded and it is quadratically upper bounded due to Assumption~\ref{ass:upper_bound} such that $\xi^s$ is exponentially stable in closed loop due to standard Lyapunov arguments~\citep{rawlings2020model}.
$\hfill\square$
\end{pf}

Theorem~\ref{thm:stability} shows that the proposed MPC scheme based on Problem~\eqref{eq:DD_MPC} exponentially stabilizes the desired setpoint $\xi^s$ and hence, the closed-loop input/output exponentially converges to $u^s$/$y^s$.
Since the cost of Problem~\eqref{eq:DD_MPC} is only positive semidefinite in the minimal state $x$, the proof relies on detectability to construct a Lyapunov function based on IOSS arguments, compare also~\cite[Theorem 2.24]{rawlings2020model}, \cite{berberich2021guarantees,grimm2005model}.
Note that Theorem~\ref{thm:stability} requires persistence of excitation of order $L+l+n$ for the input since the length of the constructed input-output trajectory is $L+l$ due to the initial conditions in~\eqref{eq:DD_MPC_init}.

\begin{remark}\label{rk:alternative_term_ing}
We note that, if the terminal set constraint in~\eqref{eq:DD_MPC_term} is dropped, i.e., $\Xi_f=\mathbb{R}^{n_{\xi}}$, closed-loop stability as shown in Theorem~\ref{thm:stability} holds locally with some region of attraction as long as the terminal cost matrix $P$ satisfies the conditions in Assumption~\ref{ass:term_ingredients}, compare~\cite{limon2006stability}.
For the special case of open-loop stable systems with no output constraints, i.e., $\mathbb{Y}=\mathbb{R}^p$, Assumption~\ref{ass:term_ingredients} can always be satisfied with $K=0$ and without the terminal set constraint in~\eqref{eq:DD_MPC_term}, i.e., $\Xi_f=\mathbb{R}^{n_{\xi}}$, compare~\cite{rawlings1993stability}.
In this case, Problem~\eqref{eq:DD_MPC} is globally feasible and Theorem~\ref{thm:stability} ensures global exponential stability.
\end{remark}

\subsection{Data-driven design of terminal ingredients}\label{subsec:MPC_design}

In this section, we show how terminal ingredients satisfying Assumption~\ref{ass:term_ingredients} can be computed using only measured data and no model knowledge, based on recent results on data-driven controller design by~\cite{berberich2020combining}.
To this end, we write~\eqref{eq:sys_xi} as a linear fractional transformation
\begin{align}\label{eq:sys_H2}
\left[\begin{array}{c}\xi_{k+1}\\\hline z_k\end{array}\right]&=
\left[\begin{array}{c|ccc}A'&B_w&B'\\\hline
\begin{bmatrix}I\\0\end{bmatrix}&0&\begin{bmatrix}0\\I\end{bmatrix}\end{array}\right]
\left[\begin{array}{c}\xi_{k}\\\hline w_k\\u_k\end{array}\right],\\\nonumber
w_k&=\Delta z_k,
\end{align}
where $w\mapsto z$ represents an additional uncertainty channel capturing all unknown elements of~\eqref{eq:sys_xi} (compare~\cite{zhou1996robust}).
More precisely, $A'$, $B'$, $B_w=\begin{bmatrix}0&I\end{bmatrix}^\top$ are known matrices according to the structure in~\eqref{eq:sys_xi_detailed}, and
\begin{align*}
\Delta=\begin{bmatrix}G_l&\dots&G_1&F_l&\dots&F_1&D\end{bmatrix}
\end{align*}
contains the unknown system parameters.
We factorize $T_y^\top QT_y=Q_r^\top Q_r$, $R=R_r^\top R_r$ with $T_y$ as in Section~\ref{subsec:prelim_extended}, and we define the data-dependent matrices
\begin{align*}
\Xi&\coloneqq\begin{bmatrix}\xi^d_l&\xi^d_{l+1}&\dots&\xi^d_{N-1}\end{bmatrix},\>
\Xi_+\coloneqq\begin{bmatrix}\xi^d_{l+1}&\xi^d_{l+2}&\dots&\xi^d_{N}\end{bmatrix},\\
U&\coloneqq\begin{bmatrix}u^d_l&u^d_{l+1}&\dots&u^d_{N-1}\end{bmatrix},\>
Z\coloneqq\begin{bmatrix}\Xi\\U\end{bmatrix},
\end{align*}
where $\{\xi^d_k\}_{k=l}^{N-1}$ is the extended state trajectory corresponding to the available input-output measurements $\{u_k^d,y_k^d\}_{k=0}^{N-1}$.
Using that $B_w^\top$ is the Moore-Penrose inverse of $B_w$, \cite[Proposition 1]{berberich2020combining} implies the following data-dependent bound on the ``uncertainty'' $\Delta$:
\begin{align}\label{eq:P_Delta_def}
\begin{bmatrix}\Delta^\top\\I\end{bmatrix}^\top
\underbrace{\begin{bmatrix}-ZZ^\top&Z\Xi_+^\top B_w\\B_w^\top\Xi_+Z^\top&-B_w^\top\Xi_+\Xi_+^\top B_w\end{bmatrix}}_{P_{\Delta}^w\coloneqq}
\begin{bmatrix}\Delta^\top\\I\end{bmatrix}
\succeq0.
\end{align}
Moreover, we abbreviate
\begin{align*}
\bar{P}_{\Delta}^w\coloneqq\begin{bmatrix}0&I\\B_w^\top&0\end{bmatrix}^\top P_{\Delta}^w\begin{bmatrix}0&I\\B_w^\top&0\end{bmatrix}.
\end{align*}
\begin{proposition}\label{prop:term_ing}
Suppose there exist $\mathcal{X}\succ0$, $\Gamma\succ0$, $M$, $\tau\geq0$, $\gamma>0$ such that $\mathrm{trace}(\Gamma)<\gamma^2$, $\begin{bmatrix}\Gamma&I\\I&\mathcal{X}\end{bmatrix}\succ0$, and
\begin{align}\label{eq:prop_term_ing_LMI}
&\begin{bmatrix}
\left(\tau \bar{P}_{\Delta}^w-\begin{bmatrix}\mathcal{X}&0\\0&0\end{bmatrix}\right)&
\begin{bmatrix}A'\mathcal{X}+B'M\\\mathcal{X}\\M\end{bmatrix}&0\\
\star&-\mathcal{X}&\begin{bmatrix}Q_r\mathcal{X}\\R_rM\end{bmatrix}^\top\\
\star&\star&-I\end{bmatrix}\prec0,
\end{align}
define $P\coloneqq \mathcal{X}^{-1}-T_y^\top QT_y$, $K\coloneqq M\mathcal{X}^{-1}$, and choose $\beta>0$ such that
\begin{align}\label{eq:prop_term_ing_region}
\Xi_f\coloneqq&\{\xi\in\mathbb{R}^{n_\xi}\mid \lVert\xi-\xi^s\rVert_P^2\leq \beta\}\subseteq\mathbb{U}^l\times\mathbb{Y}^l.
\end{align}
Then, Assumption~\ref{ass:term_ingredients} holds with $P$, $K$, and $\Xi_f$.
\end{proposition}
\begin{pf}
By applying the Schur complement to the right-lower block of~\eqref{eq:prop_term_ing_LMI} and subsequently left- and right-multiplying $\mathrm{diag}(I,\mathcal{X}^{-1})$, we obtain
\begin{align}\label{eq:prop_term_ing_proof1}
\begin{bmatrix}\tau \bar{P}_{\Delta}^w-\begin{bmatrix}\mathcal{X}&0\\0&0\end{bmatrix}&
\begin{bmatrix}A'+B'K\\I\\K\end{bmatrix}\\
\star&T_y^\top QT_y+K^\top RK-\mathcal{X}^{-1}\end{bmatrix}\prec0.
\end{align}
Applying the Schur complement to the right-lower block of~\eqref{eq:prop_term_ing_proof1} and re-arranging terms leads to $\mathcal{Q}\coloneqq \mathcal{X}^{-1}-T_y^\top QT_y-K^\top RK\succ0$ as well as
\begin{align*}
\left[\begin{array}{cc}
I&0\\(A'+B'K)^\top&\begin{bmatrix}I\\K\end{bmatrix}^\top\\\hline
0&I\\B_w^\top&0
\end{array}
\right]^\top
\left[
\begin{array}{c|c}
\begin{matrix}-\mathcal{X}&0\\0&\mathcal{Q}^{-1}\end{matrix}&0\\\hline
0&\tau P_{\Delta}^w
\end{array}
\right]
\left[
\begin{array}{cc}
\star\\
\star\\\hline
\star\\
\star
\end{array}
\right]\prec0.
\end{align*}
Using the full-block S-procedure (compare~\cite{scherer2001lpv}) and $\tilde{A}=A'+B_w\Delta\begin{bmatrix}I\\0\end{bmatrix}$, $\tilde{B}=B'+B_w\Delta\begin{bmatrix}0\\I\end{bmatrix}$, where $\Delta$ satisfies~\eqref{eq:P_Delta_def}, the latter inequality implies
\begin{align}
-\mathcal{X}+(\tilde{A}+\tilde{B}K)\mathcal{Q}^{-1}(\tilde{A}+\tilde{B}K)^\top\prec0.
\end{align}
Applying the Schur complement twice, this in turn implies
\begin{align}\label{eq:prop_term_ing_proof2}
(\tilde{A}+\tilde{B}K)^\top \mathcal{X}^{-1}(\star)
-\mathcal{X}^{-1}+T_y^\top QT_y+K^\top RK\prec0.
\end{align}
Let now $\xi^+=(\tilde{A}+\tilde{B}K)\xi$, $y=(\tilde{C}+\tilde{D}K)\xi$ for some $\xi\in\mathbb{R}^{n_\xi}$.
We then have
\begin{align*}
\lVert\xi^+\rVert_P^2
&\stackrel{\eqref{eq:prop_term_ing_proof2}}{\leq}
\lVert\xi\rVert_{\mathcal{X}^{-1}}-\lVert\xi\rVert_{T_y^\top QT_y+K^\top RK}^2
-\lVert\xi^+\rVert_{T_y^\top QT_y}^2\\
&=\lVert\xi\rVert_P-\lVert y\rVert_Q^2-\lVert\xi\rVert_{K^\top RK}^2,
\end{align*}
which proves Part (iii) of Assumption~\ref{ass:term_ingredients}.
Part (i) is a simple consequence of Part (iii).
Finally, Part (ii) follows from invariance of $\Xi_f\subseteq\mathbb{U}^l\times\mathbb{Y}^l$.
$\hfill\square$
\end{pf}

Proposition~\ref{prop:term_ing} provides a procedure to compute terminal ingredients satisfying Assumption~\ref{ass:term_ingredients}, using no model knowledge but only input-output data.
The proof relies on the full-block S-procedure~(\cite{scherer2001lpv}), which is used to robustify the model-based matrix inequality~\eqref{eq:ass_term_ingredients} against the unknown parameters $\Delta$ satisfying the known quadratic bound~\eqref{eq:P_Delta_def} derived by~\cite{berberich2020combining}, compare also~\cite{waarde2020from} for a similar approach to data-driven $\mathcal{H}_2$-state-feedback design.
The design procedure in Proposition~\ref{prop:term_ing} is an LMI feasibility problem which can be solved efficiently in practice.
Although any feasible solution $\mathcal{X}\succ0$, $M$, $\tau\geq0$ satisfying~\eqref{eq:prop_term_ing_LMI} leads to terminal ingredients satisfying Assumption~\ref{ass:term_ingredients}, the additional conditions $\mathrm{trace}(\Gamma)<\gamma^2$, $\begin{bmatrix}\Gamma&I\\I&\mathcal{X}\end{bmatrix}\succ0$ for some $\Gamma\succ0$, $\gamma>0$ lead to a smaller terminal cost.
To be precise, the optimal solution minimizing $\gamma$ subject to the conditions in Proposition~\ref{prop:term_ing} (provided it exists) is equivalent to a linear quadratic regulator for the system~\eqref{eq:sys_xi}, compare~\cite{berberich2020combining} for details.
Hence, this optimal solution provides a standard choice for terminal ingredients, compare~\cite{chen1998quasi}, leading to a good closed-loop performance.
For a reliable numerical implementation, we perform a bisection algorithm over $\gamma>0$ subject to the conditions in Proposition~\ref{prop:term_ing}.

\begin{remark}\label{rk:noise}
The proposed MPC scheme as well as the theoretical analysis in Section~\ref{subsec:MPC_stab} and the design of the terminal ingredients explained above all require that the measured data are noise-free.
However, an extension of the presented results to noisy data is straightforward, where the noise may enter both the offline data used to build the Hankel matrices in~\eqref{eq:DD_MPC_hankel} and to compute the terminal ingredients via Proposition~\ref{prop:term_ing} as well as the online data used to specify initial conditions in~\eqref{eq:DD_MPC_init}.
In particular, it is shown in~\cite{berberich2021linearpart2} that (offline and online) output measurement noise in data-driven MPC can be viewed as an input disturbance for the corresponding nominal (noise-free) data-driven MPC scheme.
Since model-based MPC with terminal ingredients possesses inherent robustness properties~\citep{yu2014inherent}, we conjecture that the closed loop is practically exponentially stable w.r.t. the noise and disturbance level, provided that Problem~\eqref{eq:DD_MPC} is slightly modified (by including a slack variable in~\eqref{eq:DD_MPC_hankel} and a regularization on $\alpha(t)$ in the cost, cf.~\cite{berberich2021linearpart2} for details).
%In particular, following similar arguments as in~\cite{berberich2021guarantees}, we conjecture that after a suitable modification of the MPC scheme (including a slack variable in~\eqref{eq:DD_MPC_hankel} and a regularization on $\alpha(t)$ in the cost), the closed loop is practically exponentially stable w.r.t. the noise and disturbance level.
As a noteworthy advantage, this is possible based on a one-step MPC scheme as in Algorithm~\ref{alg:MPC}, whereas the robustness guarantees from~\cite{berberich2021guarantees,berberich2021linearpart2} require the application of a multi-step MPC scheme due to the presence of terminal equality constraints.
%Moreover, if compared to such terminal equality constraints, MPC schemes using general terminal ingredients possess a larger closed-loop region of attraction and an improved robustness~(\cite{yu2014inherent}).
Finally, the (offline) design by~\cite{berberich2020combining} which forms the basis for Proposition~\ref{prop:term_ing} also applies in the presence of noise and hence, our design of terminal ingredients can be carried out in this case.
More precisely, replacing the definition of $P_{\Delta}^w$ in~\eqref{eq:P_Delta_def} by
\begin{align*}
P_{\Delta}^w=\begin{bmatrix}-ZZ^\top&Z\Xi_+^\top B_w\\
B_w^\top \Xi_+Z^\top&\bar{d}I-B_w^\top \Xi_+\Xi_+^\top B_w
\end{bmatrix}
\end{align*}
for some $\bar{d}>0$, Proposition~\ref{prop:term_ing} leads to terminal ingredients which satisfy Assumption~\ref{ass:term_ingredients} robustly for all systems~\eqref{eq:sys_xi} which are consistent with the measured data, when assuming that data generated offline are perturbed by process noise $\{d_k\}_{k=0}^{N-1}$ and measurement noise $\{\varepsilon_k\}_{k=0}^{N}$ satisfying a bound $\sum_{k=0}^{N-1}\lVert d_k\rVert_2^2+\lVert\varepsilon_{k+1}-\tilde{A}\varepsilon_k\rVert_2^2\leq\bar{d}$ (compare~\cite{berberich2020combining} for details).
To summarize, all results in this paper can be extended to the realistic scenario where the measured data are affected by noise and the system is subject to disturbances.
\end{remark}

\begin{remark}\label{rk:sysID}
In case of noise-free data, an \emph{indirect} approach consisting of system identification and model-based MPC is a simple alternative to our data-driven MPC scheme.
In particular, with persistently exciting data, the system matrices can be computed exactly and thus, a model-based design of terminal ingredients and implementation of MPC can be carried out analogous to Proposition~\ref{prop:term_ing} and Algorithm~\ref{alg:MPC}, respectively.
The resulting computational complexity is comparable to that of the \emph{direct} approach proposed in this paper.
On the other hand, as discussed in more detail in Remark~\ref{rk:noise}, our data-driven MPC can be extended to noisy data, in which case it may be preferable over the indirect approach due to the challenges in obtaining tight estimation bounds from data~\citep{matni2019self}.
A theoretical comparison of indirect and direct data-driven MPC is an interesting direction for future research, cf.~\cite{doerfler2021bridging,krishnan2021on} for existing open-loop results.
\end{remark}

In the remainder of this section, we investigate when a crucial necessary condition for feasibility of the conditions in Proposition~\ref{prop:term_ing} is fulfilled.
Since~\eqref{eq:prop_term_ing_LMI} is a strict LMI, its feasibility requires $\tau \bar{P}_{\Delta}^w-\begin{bmatrix}\mathcal{X}&0\\0&0\end{bmatrix}\prec0$, which implies $-ZZ^\top\prec0$.
This in turn requires that the data matrix $Z=\begin{bmatrix}\Xi^\top&U^\top\end{bmatrix}^\top$ has full row rank.
The latter condition is related to persistence of excitation (Definition~\ref{def:pe}).
In particular,~\cite[Corollary 2 (ii)]{willems2005note} implies that $Z$ has full row rank if $(\tilde{A},\tilde{B})$ is controllable and $u^d$ is persistently exciting of a sufficiently high order.
Thus, one possibility to ensure that~\eqref{eq:prop_term_ing_LMI} is feasible via a sufficiently rich input signal is to require that $(\tilde{A},\tilde{B})$ is controllable, which need not be the case in general, even if $(A,B)$ is controllable.
In order to determine when this sufficient condition for $Z$ having full row rank holds, the following result provides an equivalent condition for controllability of $(\tilde{A},\tilde{B})$ based on the observability matrix $\Phi_l$ in~\eqref{eq:obsv}.

\begin{lemma}\label{lem:extended_ctrb}
The pair $(\tilde{A},\tilde{B})$ is controllable if and only if $\Phi_l$ is square.
\end{lemma}
\begin{pf}
We show that an arbitrary $\xi_k\in\mathbb{R}^{n_\xi}$ is reachable from any initial condition if and only if $\Phi_l$ is square.
The system dynamics~\eqref{eq:sys} imply for any $k\geq n_{\xi}$
\begin{align}\label{eq:lem_ctrb_proof1}
\xi_k=&\begin{bmatrix}u_{[k-l,k-1]}\\\Phi_lx_{k-l}+\Gamma_l u_{[k-l,k-1]}\end{bmatrix}
\end{align}
with a suitably defined matrix $\Gamma_l$.
%
%\begin{align*}
%\Gamma_l&\coloneqq\begin{bmatrix}D&0&\dots&0\\
%CB&D&\ddots&\vdots\\
%\vdots&\ddots&\ddots&0\\
%CA^{l-1}B&\dots&CB&D
%\end{bmatrix}.
%\end{align*}
%
By observability, i.e., full column rank of $\Phi_l$, we have $pl\geq p\underline{l}\geq n$ and hence, $k-l\geq n_{\xi}-l\geq n$.
Thus, using controllability of $(A,B)$, the system can be steered to an arbitrary state $x_{k-l}$ from any initial condition $x_0$ by appropriately selecting the input $u_{[0,k-l-1]}$.
Therefore, the state $\xi_k$ is reachable if and only if there exist an input $u_{[k-l,k-1]}$ and a state $x_{k-l}$ such that~\eqref{eq:lem_ctrb_proof1} holds.
The first block row of~\eqref{eq:lem_ctrb_proof1} is trivially satisfied by definition of $\xi_k$.
Using $u_{[k-l,k-1]}=\begin{bmatrix}I&0\end{bmatrix}\xi_k$, the second block row of~\eqref{eq:lem_ctrb_proof1} takes the form
\begin{align}\label{eq:lem_ctrb_proof2}
\begin{bmatrix}-\Gamma_l&I\end{bmatrix}\xi_k=\Phi_lx_{k-l}.
\end{align}
Clearly, there exists a state $x_{k-l}$ for any given $\xi_k$ such that~\eqref{eq:lem_ctrb_proof2} holds if and only if $\begin{bmatrix}-\Gamma_l&I\end{bmatrix}\xi_k$ lies in the image of $\Phi_l$ for any $\xi_k\in\mathbb{R}^{n_\xi}$.
Since $\begin{bmatrix}-\Gamma_l&I\end{bmatrix}$ has full row rank, this means that $\Phi_l$ has full row rank.
Since $\Phi_l$ has full column rank by the assumption that $(A,C)$ is observable, full row rank of $\Phi_l$ is equivalent to $\Phi_l$ being square.
$\hfill\square$
\end{pf}

Lemma~\ref{lem:extended_ctrb} shows that, under the present assumptions, the extended dynamics~\eqref{eq:sys_xi} are controllable if and only if the matrix $\Phi_l$ is square.
This means that full row rank of $Z$, which is a necessary condition for the presented design of terminal ingredients and thus the proposed MPC scheme, can be guaranteed if $\Phi_l$ is square.
The latter condition requires $l=\underline{l}$, i.e., the lag of the system is known and chosen for defining the extended state $\xi$, as well as $pl=n$.
For single-output systems with $p=1$, this in turn requires that $l=n$, i.e., the system order is known exactly.
For multiple-output systems on the other hand, the condition $pl=n$ does in general not hold even if both the system order and the lag are known exactly.
Addressing this issue and designing data-driven controllers for the extended system~\eqref{eq:sys_xi} in case it is not controllable is an interesting issue for future research, independently of its application to data-driven MPC in this paper.

\section{Numerical example}\label{sec:example}
\begin{figure}
		\begin{center}
		\subfigure
		{\includegraphics[width=0.49\textwidth]{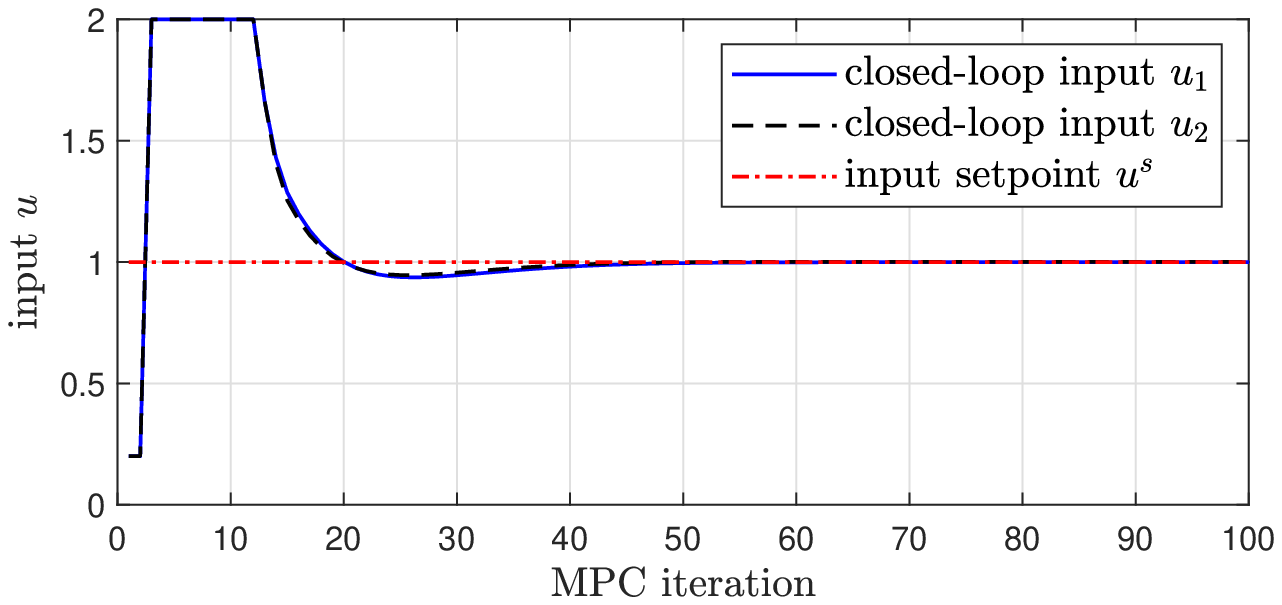}}
		\subfigure
		{\includegraphics[width=0.49\textwidth]{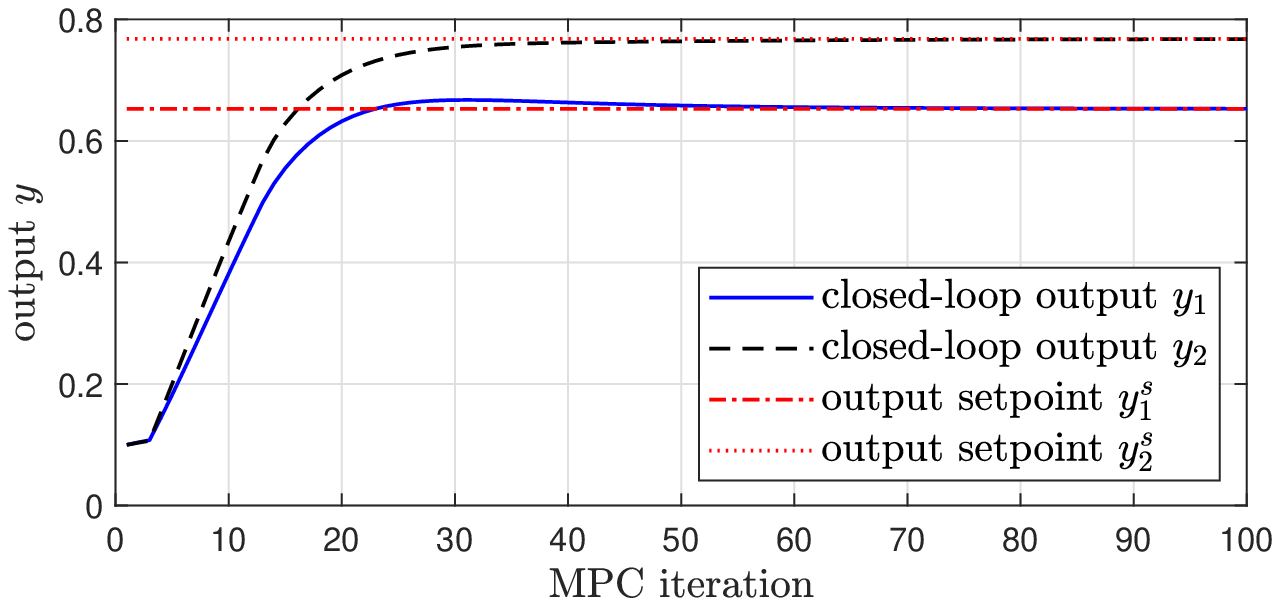}}
		\end{center}
		\caption{Closed loop under the data-driven MPC scheme.}	\label{fig:four_tank}
\end{figure}

We apply our MPC scheme to a linearized version of the four-tank system considered by~\cite{raff2006nonlinear}, see~\cite{berberich2021guarantees} for the numerical values of $A$, $B$, $C$, $D$.
%
%We revisit the example considered by~\cite{berberich2021guarantees} and apply our MPC scheme to System~\eqref{eq:sys} with
%%
%\begin{align*}
%A&=\begin{bmatrix}0.921&0&0.041&0\\
%0&0.918&0&0.033\\
%0&0&0.924&0\\
%0&0&0&0.937\end{bmatrix},\>\>
%B=\begin{bmatrix}0.017&0.001\\
%0.001&0.023\\
%0&0.061\\
%0.072&0\end{bmatrix},\\
%C&=\begin{bmatrix}1&0&0&0\\
%0&1&0&0\end{bmatrix},\>\>D=0.
%\end{align*}
%%
%These dynamics correspond to a linearized version of the four-tank system considered by~\cite{raff2006nonlinear}.
We consider the prediction horizon $L=15$, cost matrices $Q=I$, $R=5\cdot 10^{-3}I$, an input constraint set $\mathbb{U}=[-2,2]^2$, but no output constraints, i.e., $\mathbb{Y}=\mathbb{R}^p$, and the input-output setpoint $(u^s,y^s)=\left(\begin{bmatrix}1\\1\end{bmatrix},\begin{bmatrix}0.65\\0.77\end{bmatrix}\right)$.
Moreover, we regularize $\alpha(t)$ in the cost of Problem~\eqref{eq:DD_MPC} by penalizing $10^{-4}\cdot\lVert\alpha(t)\rVert_2^2$ for better numerical stability and robustness (compare~\cite{coulson2020distributionally,berberich2021guarantees}).
We generate data of length $N=100$ of the above system by applying an input sampled uniformly from $u_k\in[-1,1]^2$.
First, we note that the observability matrix $\Phi_l$ of this system is square, i.e., the extended system is controllable (compare Lemma~\ref{lem:extended_ctrb}), for $l=\underline{l}=2$, i.e., if $l$ is chosen as the lag of the system.
Based on this knowledge, we apply Proposition~\ref{prop:term_ing} to compute terminal ingredients satisfying Assumption~\ref{ass:term_ingredients}.
To be precise, using Yalmip~(\cite{lofberg2004yalmip}) with the solver MOSEK~(\cite{MOSEK15}), we obtain a feasible solution of the conditions in Proposition~\ref{prop:term_ing} with $K=0$ and $\gamma=52$.
We then apply the MPC scheme based on Problem~\eqref{eq:DD_MPC} with $\Xi_f=\mathbb{R}^{n_\xi}$.
Since the above system is open-loop stable, the closed loop under this MPC scheme is globally exponentially stable (compare Remark~\ref{rk:alternative_term_ing}).
The closed-loop input and output trajectories are displayed in Figure~\ref{fig:four_tank}.
It can be seen that the MPC scheme tracks the desired setpoint reliably.
On the other hand, an MPC scheme using no terminal ingredients (i.e., setting $P=0$ in~\eqref{eq:DD_MPC}) leads to an \emph{unstable} closed loop since the prediction horizon is chosen too short.
Moreover, an MPC scheme based on terminal equality constraints as proposed by~\cite{berberich2021guarantees} is not initially feasible for the initial condition $x_0=\begin{bmatrix}0.1&0.1&0.2&0.2\end{bmatrix}^\top$ with any prediction horizon $L\leq24$ due to the input constraints.
This means that $x_0$ does not lie in the (bounded) region of attraction of the scheme by~\cite{berberich2021guarantees}, in contrast to the global stability guarantees of the MPC scheme designed above.
Finally, we note that it is also possible to design a non-zero $K$ based on Proposition~\ref{prop:term_ing} resulting in a less conservative terminal penalty.
However, in this case, non-trivial terminal set constraints need to be included which reduces the region of attraction.
To conclude, our results show that enhancing data-driven MPC with non-trivial terminal ingredients improves both the theoretical properties as well as the practical performance.

\section{Conclusion}\label{sec:conclusion}
We presented a data-driven MPC scheme with terminal ingredients.
The scheme uses only measured input-output data for prediction and thus, the terminal ingredients involve an extended state consisting of consecutive input-output values.
We proved closed-loop stability and showed that the required terminal ingredients can be computed based on data without explicit model knowledge.
If compared to existing data-driven MPC schemes with closed-loop guarantees~(\cite{berberich2021guarantees}), the presented approach has superior closed-loop robustness and performance properties both in theory and practice.

\bibliographystyle{ifacconf}   
\bibliography{Literature}  

\end{document}